**May 25, 2007**

# Equiangular Surfaces, Self-Similar Surfaces, and the Geometry of Seashells

**Khristo N. Boyadzhiev**
Department of Mathematics and Statistics,
Ohio Northern University, Ada, OH 45810, USA
[k-boyadzhiev@onu.edu](k-boyadzhiev@onu.edu)

One of the most fascinating curves in the plane is the logarithmic spiral shown in Figure 1, also known as the equiangular spiral. Among its interesting properties is the fact that it is self-similar (we explain this in a moment). In this article, we explore three-dimensional versions of these two properties: surfaces that are equiangular and those that are self-similar. It will be seen that the two concepts are independent in three dimensions but not in two. Among equiangular surfaces, there are some that are self-similar, but also some that are not. Furthermore, there are self-similar surfaces that are not equiangular. We will give examples of all three types.

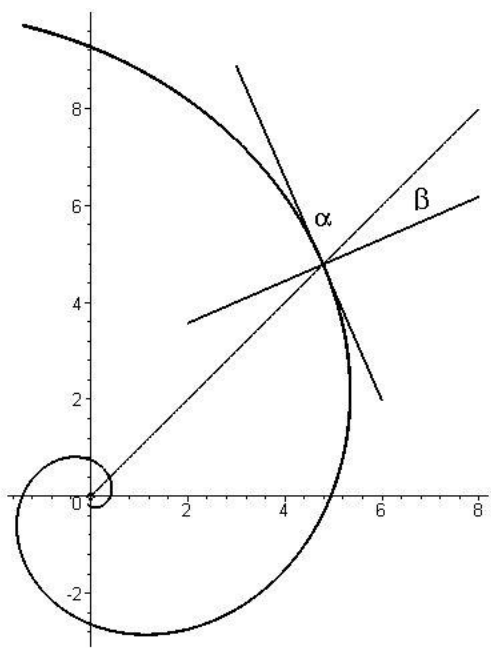

**Figure 1.** A logarithmic spiral.

A *logarithmic spiral* is defined by the polar equation

$$r = r_0 e^{k\theta} \qquad \left(-\infty < \theta < \infty\right), \tag{1}$$

where $\theta$ is the polar angle and $r = r(\theta)$ is the polar radius. Figure 1 shows an example of such a curve. This spiral is equiangular - at each point the angle $\alpha$ between the radius vector and the tangent vector is always the same, with $\cot \alpha = k$ (see Figure 1). In [2], we showed that this property leads to (1), that is, it characterizes logarithmic spirals. More about these curves can be found in [2], [3], and [9], and the references therein.

A wonderful fact about these spirals is that they appear in nature, including in shells and horns (see for example, [1], [3], [5], [8]). In particular, a cross section of the popular *nautilus pompilius* reveals a logarithmic spiral [3], [8]. This phenomenon is due to the property of self-similarity: any change of the parameter $\theta \to \theta + \gamma$ leads to a proportional curve $r_\gamma(\theta) = (r_0 e^{k\gamma}) e^{k\theta}$. This property too characterizes logarithmic spirals since the functional equation $f(s+t) = f(s)f(t)$ for one continuous function implies that the function $f$ is an exponential, that is, $f(t) = e^{at}$.

Shells, on their part, grow as self-similar structures. According to D'Arcy Thompson [8, p. 757], "In the growth of a shell, we can conceive no simpler law than this, namely, that it shall widen and lengthen in the same unvarying proportions: and this simplest of laws is that which Nature tends to follow. The shell, like the creature within it, grows in size *but does not change its shape*; and the existence of this constant relativity of growth, or constant similarity of form, is of the essence, and may be made the basis of a definition, of the equiangular spiral."

Since shells are 3-dimensional, it is interesting to investigate surfaces in 3-space that are analogues of the logarithmic spiral. Several successful constructions have been made of surfaces modeling shells (see, for instance, [4], [6], [7]). Most of these constructions are based on an equiangular spiral on a cone, called a conchospiral [2] (or helico-spiral). By moving a closed generating curve along a conchospiral, one can get a variety of shell-like surfaces, and Meinhardt [7] gives some interesting computer-generated images of seashells. However, he does not discuss equiangular properties of these surfaces.

We take a different approach, looking directly for surfaces having the desirable properties of the spiral: equiangularity and self-similarity. We say that a surface is *self-similar* if its polar equation is of the form

$$\rho(\theta, \psi) = \rho_0 e^{\mu\theta + b\psi}, \tag{2}$$

where $\rho$ is the polar radius, $\theta$ and $\psi$ are respectively the latitude and longitude, and $\mu$ and $b$ are constants. We also include all planes and cones even though they cannot be described by a polar equation of the form (2).

Clearly, replacing $\theta$ by $\theta + \Delta\theta$ or $\psi$ by $\psi + \Delta\psi$ in (2) yields a proportional shape, and so we are justified in calling such a surface self-similar. So, a natural question is this: Are such surfaces equiangular? In the next section we define and investigate equiangular surfaces, after which we show that the two types of surface are quite different.

**Equiangular Surfaces**

Based on the equiangular property of logarithmic spirals, we say that a surface is *equiangular* if it is smooth and the angle $\beta$ between the normal vector $N$ and the radius vector $X$ is the same at all points on the surface. Specific values of $\beta$ (called the *characteristic angle* of the surface) provide some degenerate examples, all of which are also self-similar:

- Spheres centered at the origin $(\beta = 0)$;

- Cones with vertex at the origin and planes through the origin $(\beta = \frac{\pi}{2})$.

More interesting examples are obtained by rotating a logarithmic spiral about a coordinate axis. For example, if the spiral (in the $xz$-plane)

$$\begin{cases} x = r_0 e^{kt} \cos t \\ z = r_0 e^{kt} \sin t \end{cases}$$

is rotated about the $z$-axis, we obtain the surface

$$\begin{cases} x = r_0 e^{kt} \cos t \cos \theta \\ y = r_0 e^{kt} \cos t \sin \theta \\ z = r_0 e^{kt} \sin t \end{cases} . \tag{3}$$

Here, $\theta$ is the angle of rotation, with $0 \le \theta \le 2\pi$, while $-\infty < t < +\infty$. Figure 2 shows the portion of this surface for $-\frac{\pi}{2} \le t \le \frac{3\pi}{2}$ and a cross-section of this portion. The entire surface consists of infinitely many nested egg-like shells. Letting $t = \psi$, we see that the polar equation of this surface is

$$\rho(\theta, \psi) = r_0 e^{k\psi}, \tag{4}$$

and therefore it is self-similar. We have already seen that it is equiangular, so it is an example of a surface with both key properties.

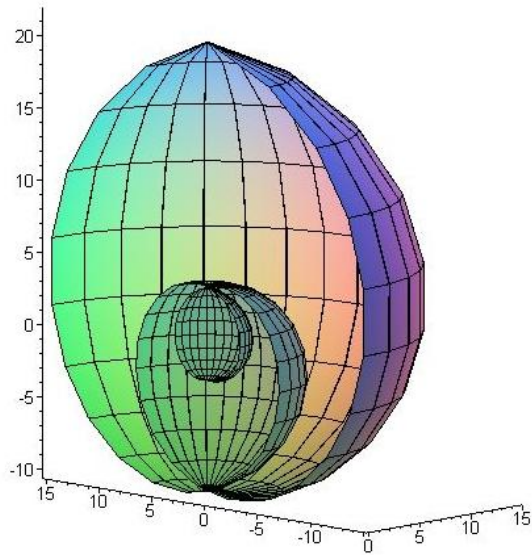

**Figure 2.** The equiangular and self-similar matryoshka-like surface of revolution.

The general equiangular surface satisfies a differential equation that will be useful in our investigation. We find it convenient to work in polar coordinates (as in our last example), with the standard relationship with rectangular coordinates being

$$\begin{cases} x = \rho \cos \psi \cos \theta \\ y = \rho \cos \psi \sin \theta \\ z = \rho \sin \psi \end{cases} . \qquad (5)$$

Where we allow $-\infty < \theta < \infty$ and $-\infty < \psi < \infty$. We adopt the standard convention of using subscripts to denote partial derivatives.

**Theorem.** *If S is an equiangular surface with polar equation $\rho = \rho(\theta, \psi)$ and with characteristic angle $\beta \neq \frac{\pi}{2}$, then*

$$\rho_\theta^2 + \rho_\psi^2 \cos^2 \psi = \rho^2 \cos^2 \psi \, \tan^2 \beta. \qquad (6)$$

Note that of the equiangular surfaces we have considered thus far, spheres and rotational spiral surfaces are included here but planes and cones are not. For spheres centered at the origin, $\beta = 0$ and $\rho_\theta = \rho_\psi = 0$. For rotational spiral surfaces, $\rho_\theta = 0$ and so $\rho = \rho_0 e^{(\pm \tan \beta)\psi}$, as in (4).

Before proving the theorem, we show that it answers the question "Which self-similar surfaces are equiangular?" posed earlier. Substituting the equation for a self-similar surface (2) into the differential equation (6) and then doing a bit of algebra, we find that $\mu^2 = (\tan^2 \beta - b^2)(\cos^2 \psi)$. Since $\psi$ varies but $\mu, \beta$, and $b$ are constant, it follows that $\mu = 0$ and $b = \pm \tan \beta$. Consequently, we have the following result.

**Corollary**. *The only self-similar surfaces that are equiangular with characteristic angle $\beta$ have polar equation $\rho(\theta, \psi) = \rho_0 e^{(\pm \tan \beta)\psi}$.*

Figure 3 presents an example of a surface which is self-similar, but not equiangular. It is defined by equation (2) with $\mu = 0.5$ and $b = 0.2$, where $0 < \theta < 2\pi$ and $-\dfrac{\pi}{2} < \psi < \dfrac{\pi}{2}$ .

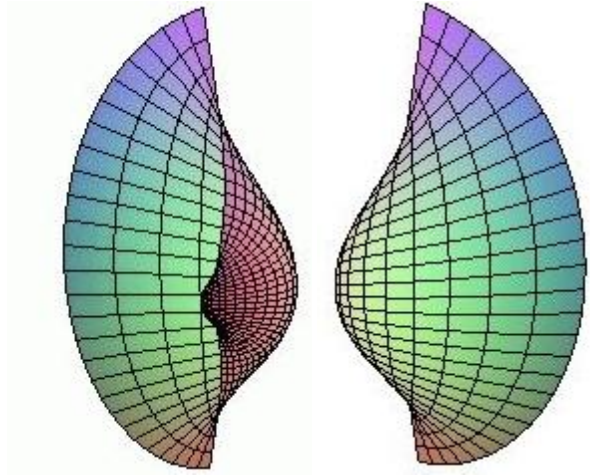

**Figure 3.** A self-similar shell which is not equiangular.

*Proof of theorem.* Let $\boldsymbol{X}(\theta, \psi) = \langle x, y, z \rangle$ be the radius vector of the surface $S$, where $x$, $y$, and $z$ are defined by (5). The normal vector $\boldsymbol{N} = \langle n_1, n_2, n_3 \rangle$ equals $\boldsymbol{X}_\theta \times \boldsymbol{X}_\psi$, so by simple computation we find

$$\begin{cases} n_1 = \rho_\psi \rho \cos\theta \cos\psi \sin\psi + \rho_\theta \rho \sin\theta + \rho^2 \cos\theta \cos^2\psi \\ n_2 = \rho_\psi \rho \sin\theta \cos\psi \sin\psi - \rho_\theta \rho \cos\theta + \rho^2 \sin\theta \cos^2\psi \\ n_3 = -\rho_\psi \rho \cos^2\psi + \rho^2 \cos\psi \sin\psi. \end{cases}$$

It follows that $\boldsymbol{N} \cdot \boldsymbol{X} = \rho^3 \cos\psi$, that $|\boldsymbol{X}|^2 = \rho^2$, and that $|\boldsymbol{N}|^2 = \rho^2(\rho_\psi^2 \cos^2\psi + \rho^2 \cos^2\psi + \rho_\theta^2)$. The property of equiangularity requires that $\dfrac{\boldsymbol{N} \cdot \boldsymbol{X}}{|\boldsymbol{N}|\,|\boldsymbol{X}|} = \cos\beta,$ which can be written as

$$\frac{\rho^2 \cos^2\psi}{\rho_\psi^2 \cos^2\psi + \rho^2 \cos^2\psi + \rho_\theta^2} = \cos^2\beta,$$

and from this we get (6) by taking reciprocals. ∎

(A natural restriction for the above computations is $\cos\psi \neq 0$. As we shall see later, for the proper solutions of (6) we have $\cos\psi > 0$ anyway.)

**Solutions of the equiangular differential equation**

Equation (6) is a first-order non-linear partial differential equation, and we now proceed to solve it. The solution will provide us with some interesting equiangular surfaces that are not self-similar.

Our solution begins with separation of variables:

$$\left(\frac{\rho_\theta}{\rho}\right)^2 = \left(\cos^2\psi\right)\left(\tan^2\beta - \left(\frac{\rho_\psi}{\rho}\right)^2\right). \tag{7}$$

We assume that our solution has the form $\rho(\theta, \psi) = u(\theta)v(\psi)$. Then (7) becomes

$$\left(\frac{u'}{u}\right)^2 = \left(\cos^2\psi\right)\left(\tan^2\beta - \left(\frac{v'}{v}\right)^2\right).$$

Since the left side depends only on θ and the right side only on $\psi$, they must be constant, say $\mu^2$. Consequently, $\frac{u'}{u} = \pm\mu$ and so

$$u(\theta) = u_0 e^{\pm\mu\theta}. \tag{8}$$

Additionally, we find that

$$\left(\frac{v'}{v}\right)^2 = \tan^2\beta - \frac{\mu^2}{\cos^2\psi}, \tag{9}$$

from which it follows that $\mu^2 \leq \tan^2\beta$. We further deduce that $\frac{v'}{v} = \pm\sqrt{\tan^2\beta - \mu^2 - \mu^2\tan^2\psi}$.

If $\mu = 0$, then $\frac{v'}{v} = \pm\tan\beta$, so $v(\psi) = v_0 e^{\pm(\tan\beta)\psi}$, and hence the surface is one of the rotational surfaces (4).

Therefore we assume that $\mu \neq 0$, and for the sake of simplicity we consider only the case $\mu > 0$. We let

$$a = \sqrt{\left(\frac{\tan \beta}{\mu}\right)^2 - 1}\,, \tag{10}$$

so that

$$\frac{v'}{v} = \pm \mu \sqrt{a^2 - \tan^2 \psi}\,. \tag{11}$$

Note that our choice of $a$ implies that $|\tan \psi| \leq a$, and hence our solution is restricted to the interval $I_a = [-\arctan a, \arctan a]$. Now, for $\psi \in I_a$, we define $h(\psi) = \mu \int_0^\psi \sqrt{a^2 - \tan^2 \lambda}\; d\lambda$. This integral can be evaluated using tables or a computer algebra system:

$$h(\psi) = \mu \left[ \sqrt{a^2 + 1}\, \arctan\left(\frac{\sqrt{a^2 + 1}\, \tan \psi}{\sqrt{a^2 - \tan^2 \psi}}\right) - \arctan\left(\frac{\tan \psi}{\sqrt{a^2 - \tan^2 \psi}}\right) \right].$$

It follows from our definition of $h(\psi)$ that $\dfrac{v'}{v} = \pm h'$, whose solution is $v(\psi) = v_0 e^{\pm h(\psi)}$ for $\psi \in I_a$. We also set $h(\psi) = \psi \tan \beta$ when $\mu = 0$. Thus we have obtained the general solution of (6) as

$$\rho(\theta, \psi) = \rho_0 e^{\pm \mu \theta \pm h(\psi)}\,, \tag{12}$$

where $-\infty < \theta < \infty$, $-\arctan a \leq \psi \leq \arctan a$, and $\mu$ and $\rho_0$ are arbitrary constants with $\rho_0 \geq 0$. Remind that $a$, as defined in (10), depends on $\beta$ and $\mu$. When $\mu = 0$ we allow $-\infty < \psi < \infty$.

Note that the "level curves" $\psi = C$ of (12) are conchospirals [2], as follows from (5).

Any solution of (12) with $\mu \neq 0$ provides an example of an equiangular surface that is not self-similar. For example, the surface for the case $\rho_0 = 1$, $\mu = 0.1$, and $a = 2$, drawn with Maple, is shown in Figure 4, with $\theta$ restricted to the interval $[0, 4\pi]$. Changing only $a$, from 2 to 20, we obtain the prolonged shape shown in Figure 5, which looks like a biological form, while changing only $\mu$, from 0.1 to 1 and keeping $a = 2$, leads to the interesting wing-like surface shown in Figure 6.

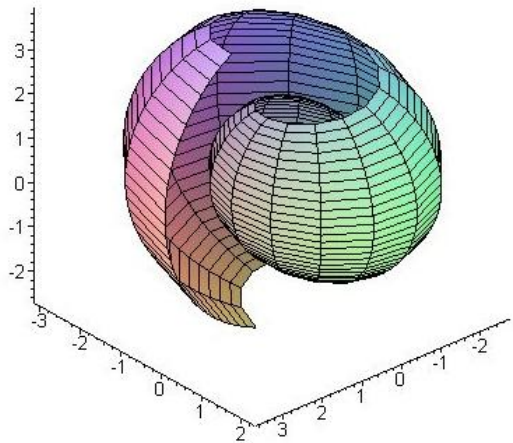

**Figure 4.** Equiangular surface, "close" to self-similar.

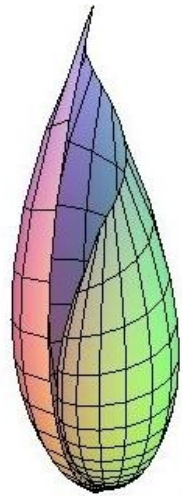

**Figure 5.** Equiangular "bud-like" surface.

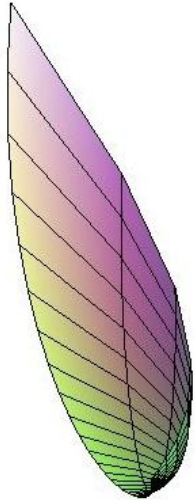

**Figure 6.** Equiangular wing.

Observe that the Maclaurin series for $h(\psi)$ is

$$h(\psi) = (\mu\,a)\psi - \frac{\mu}{6a}\,\psi^3 - \left(\frac{\mu}{a}\right)\!\left(\frac{1}{15} + \frac{1}{40a^2}\right)\psi^5 + \cdots,$$

where $\mu\,a = \sqrt{\tan^2\beta - \mu^2}$ and $\dfrac{\mu}{a} = \dfrac{\mu^2}{\sqrt{\tan^2\beta - \mu^2}}$. For fixed $0 < \beta < \dfrac{\pi}{2}$ and small $\mu$, we have

$\mu\,a \approx \tan\beta$, and so $h(\psi) \approx (\tan\beta)\psi$ and $\rho \approx \rho_0 e^{\pm\mu\theta \pm (\tan\beta)\psi}$. That is, when $\mu \neq 0$ is small, the equiangular surface (12) is close to being self-similar.

*Acknowledgment.* I would like to thank the referees for their helpful suggestions. I am also very grateful to the editor for his work on improving the exposition.